\documentclass[10pt,conference]{IEEEtran}
\usepackage{latexsym,todonotes}
\usepackage{graphicx}
\usepackage{float}
\usepackage{subfigure}
\usepackage{multirow}
\usepackage{amsmath, amssymb}
\usepackage{cite}
\usepackage{txfonts}
\usepackage{mathrsfs}
\usepackage{exscale}%
\usepackage{relsize}%
\usepackage{xcolor}
\usepackage{algorithm,algorithmicx}
\usepackage{algpseudocode} 
\usepackage{tikz,mathpazo,pgfplots}
\usepackage{cite}
\usepackage{flowchart}
\usetikzlibrary{arrows,calc,positioning}
\usepackage{tabularx}%
\usepackage{textcomp}
\usepackage{booktabs}%
\usepackage[switch]{lineno}%
\usepackage[bookmarks=false]{hyperref}

\begin{document}
\title{Sequence Q-Learning Algorithm for Optimal Mobility-Aware User Association \\
\thanks{This work is partly supported by UIC Research Grant R201911 and Zhuhai Basic and Applied Basic Research Foundation Grant ZH22017003200018PWC.}
}

\author{\IEEEauthorblockN{ Wanjun Ning\IEEEauthorrefmark{1}, Zimu Xu\IEEEauthorrefmark{1}, Jingjin Wu\IEEEauthorrefmark{1},  Tiejun Tong\IEEEauthorrefmark{2}}
\IEEEauthorblockA{\IEEEauthorrefmark{1}Department of Statistics, BNU-HKBU United International College, Zhuhai, Guangdong, P. R. China \\
\IEEEauthorrefmark{2}Department of Mathematics, Hong Kong Baptist University, Kowloon, Hong Kong SAR, P. R. China }
Email: q030202601@mail.uic.edu.cn; r130202603@mail.uic.edu.cn; jj.wu@ieee.org; tongt@hkbu.edu.hk\\}

\maketitle
\begin{abstract}
We consider a wireless network scenario applicable to metropolitan areas with developed public transport networks and high commute demands, where the mobile user equipments (UEs) move along fixed and predetermined trajectories and request to associate with millimeter-wave (mmWave) base stations (BSs). An effective and efficient algorithm, called the Sequence Q-learning Algorithm (SQA), is proposed to maximize the long-run average transmission rate of the network, which is an NP-hard problem.  Furthermore, the SQA tackles the complexity issue by only allowing possible re-associations (handover of a UE from one BS to another) at a discrete set of decision epochs and has polynomial  time complexity. This feature of the SQA also restricts too frequent handovers, which are considered highly undesirable in mmWave networks. Moreover, we demonstrate by extensive numerical results that the SQA can significantly outperform the benchmark algorithms proposed in existing research by taking all UEs’ future trajectories and possible decisions into account at every decision epoch. 

\end{abstract}

\section{Introduction}
In next-generation wireless networks, access points such as mobile base stations (BSs) are densely deployed to provide redundant coverage and fault tolerance if some of the BSs fail. As a result, mobile user equipments (UEs) can choose from multiple BSs for accessing the core network. Studies on mobility-aware user association, which focus on choosing the proper BS for each UE to improve network performance as reflected by certain Quality of Service (QoS) metrics, have thus become a popular research topic~\cite{liu2016user,tabassum2019fundamentals}. 

Operating on frequency bands between 30 and 300 GHz, mmWave communications can provide much higher transmission rates than communications on sub-6GHz bands~\cite{andrews2016modeling, Dong2020}. However, the susceptibility to physical obstacles for transmissions on mmWave bands creates new issues for research on user association strategies~\cite{andrews2016modeling,choi2020mobility}. Meanwhile, a typical mmWave BS covers a much smaller area than a conventional sub-6GHz BS, leading to more frequent re-associations (handovers between different BSs)~\cite{cacciapuoti2017mobility}. As handovers incur extra power consumption and are more likely to result in connection failures, studies on user association strategies in a mobility-aware mmWave environment must account for the number of handovers in addition to the traditional objectives and constraints.

Most existing studies on mobility-aware user association consider that UEs move along random trajectories following certain statistical distributions~\cite{arshad2017mobility,choi2020mobility,sun2017smart}. However, as most actual UE movements cannot fit in theoretical distributions~\cite{shayea2020key}, the practical value of strategies under this assumption is limited. 

This paper focuses on a scenario applicable for cities with developed public transport networks and high commute demands. Statistics show that more than 45\% of people choose public transport for commuting to and from work in metropolises such as Singapore, London, and Hong Kong~\cite{dixon20192019}. For mobile users taking subways or public buses, it is possible to estimate their UEs' future movements and positions rather accurately for a relatively long time. Due to the large proportion of such users and their usage demands, it is reasonable to propose a user association strategy for the UEs that  move along a predetermined trajectory.

The future trajectories of such UEs could be exploited by user association strategies to improve the performance further. More specifically, we will only allow a UE to handover from its current serving BS to another BS at specific points on its trajectory when certain triggering conditions are met. This approach could avoid frequent handovers and simultaneously reduce the complexity of the user association algorithm, as association decisions only need to be made at a discrete set of \emph{decision epochs} when the triggering conditions are potentially fulfilled. On the other hand, given that a BS’s available bandwidth is limited, it is beneficial to take other UEs’ possible future actions based on their trajectories into account to improve the long-run performance.

Our proposed approach to exploit the future information is applying the reinforcement learning (RL) technique. RL is a machine learning technique that learns system information from the interaction between agents and the environment to solve decision-making problems such as user association problems. For example, a deep RL (DRL)-based algorithm, Deep Q-Network, was demonstrated in~\cite{zhao2019deep} to make user association decisions regardless of mobility of UEs. A distributed method called multi-agent DQN with recurrent neural networks was proposed in~\cite{sana2020multi} to optimize user association decisions, considering channel dynamics and changing rate demand of UEs. Guo \textit{et al.}~\cite{guo2020joint} propose a multi-agent proximal policy optimization to solve the joint optimization problem of handover control and power allocation in two-layer heterogeneous cellular networks and show it is effective in small-scale experiments. Sun \textit{et al.}~\cite{sun2017smart} demonstrated a handover strategy based on RL called SMART, considering mmWave channel characteristics and UE's QoS requirements. A Q-Learning-based handover strategy was proposed in~\cite{khosravi2020learning}, which maximized the trajectory rate but ignored the interactions between UE decisions. 

The new algorithm proposed in this work is called the Sequence Q-learning Algorithm (SQA). It aims at maximizing the long-run average transmission rate of the network, a commonly used QoS metric. In particular, the SQA considers interactions between different UEs based on the predicted trajectories, which traditional Q-learning could not account due to the curse of dimensionality of the state space. Besides, the SQA can efficiently make use of future information and achieve significant performance improvement in moderate scale experiments, in which DRL approaches (e.g.~\cite{guo2020joint}) require extremely massive training steps and hardly converge. 

The main contribution of this paper is to propose an efficient SQA for mobility-aware user association, aiming at optimizing the long-run average transmission rate of the network while limiting the number of re-association to a relatively low level. The SQA is specifically designed for the scenario where the trajectories of UEs are known or can be reasonably predicted by the mobile service operators. By taking advantage of the future movement of UEs, the SQA explores appropriate weights for each possible association action based on the classical Q-learning. The weights reflect the possible impact of future movements and decisions of all UEs and are included in the adjusted action-value functions to identify the optimal BS associated with a moving UE at each decision epoch. We will demonstrate that, compared with 
state-of-the-art benchmark algorithms, the SQA achieves much better performance in terms of the long-run average transmission rate of the network.

The remainder of this paper is organized as follows. Section~\ref{sec:model} describes the system model and formulates the optimization problem. Section~\ref{sec:algorithm} presents the proposed SQA in detail. We describe the experimental setup and present the numerical results in Section~\ref{sec:experiment}, where our proposed algorithm and four state-of-the-art benchmark algorithms are compared. Section~\ref{sec:conclusion} concludes the paper.

\section{System Model and Problem Formulation}~\label{sec:model}
Let $ \mathscr{ M }= \left\lbrace 1,…,M\right\rbrace  $ denote the set of BSs, and 
$ \mathscr{ N } =\left\lbrace 1,…,N\right\rbrace $ represent the set of UEs. Accordingly, there are in total $M$ BSs and $N$ UEs randomly and uniformly distributed in the system.
 We denote $dis_{m,n}(t)$, where $m \in \mathscr{M}$ and $n \in \mathscr{N}$, as the Euclidean distance between BS $m$ and UE $n$ at time $t$. A UE is associated to one BS at any time. 

Consider $\mathscr{D}(n,t)$ as the set of candidate BSs available to be associated by UE $n$ at time $t$. As frequent handover of UEs is considered undesirable  for mobile networks~\cite{cacciapuoti2017mobility}, the re-association would only occur at decision epochs when certain triggering conditions are met. Considering the noise-limited nature of mmWave network~\cite{elshaer2016downlink}, we associate the triggering conditions with the Signal-to-Noise Ratio (SNR). Specifically, for a UE, decision epochs are the moments when the SNR from its currently serving BS falls lower than a certain threshold, or when the SNR from a neighboring BS 
becomes higher than a certain threshold~\cite{3GPPtrigger}. Assuming that such changes do not occur concurrently for more than one UE, $ \mathscr{T} =\left\lbrace t_1, t_2, …\right\rbrace$ is denoted as the set of decision epochs, where each $t_i \in \mathscr{R}$ is a specific time point when a change in $\mathscr{D}(n,t_i)$ occurs for a UE $n \in \mathscr{N}$. We assume $t_i < t_j$ for all $i < j$, thus a decision epoch $t_i$ uniquely identifies the UE $n$ that experiences a change in $\mathscr{D}(n,t_i)$. We denote this relationship as $n = U(t_i)$. For notational simplicity, we define $D(t_i) := \mathscr{D}(n,t_i)$ as the decision set at decision epoch $t_i$. 

As the previous section explains, we consider a scenario where all UEs move along predetermined trajectories. 
Accordingly, when user association decisions are made at $t_i$, the future locations of UEs at following times can be taken into account. We further assume that the latency required for a UE to obtain the states of all BSs in the network is negligible compared to the time gap between any two consecutive decision epochs, such that the network state may only change at decision epochs.

Two important evaluation metrics in optimization problems in user association are the long-run average transmission rate and the number of handovers. 
While the number of handovers is controlled by only allowing potential handovers at decision epochs, we now focus on maximizing the long-run average rate of transmission.

To formulate the optimization problem, we define another function, $\text{last}(t)$, that returns the latest decision epoch $t_i \in \mathscr{T}$ before time $t$. Notably, $\text{last}(t_i)$, where $t_i$ itself is a decision epoch, will return the last decision epoch $t_{i-1}$. Recall that $U(t)$ returns the corresponding UE $n \in \mathscr{N}$ whose set of candidate BSs is changed at $t$. 
In this way, for any $t$, we can identify a unique BS, $A_i \in D(t_i)$, where $t_i = \text{last}(t)$, that UE $U(t)$ decides to associate with at $t_i$. We further denote $\text{cnt}_m(t)$ as the number of UEs associated with BS $m$ at $t$.

We further assume that the bandwidth of a BS is evenly distributed to all associating UEs. The transmission rate of UE $n$ associated with BS $m$ at time $t$ is
\begin{equation}
	r_{m,n}(t)= \frac{B_m}{\text{cnt}_m(t)}\log_2\left(1+\frac{P_{m,n} dis_{m,n}(t)^{-\beta}}{N_0}\right),
	\label{eqn:transmission_rate}
\end{equation} 
where $P_{m,n}$ is the transmission power for BS $m$ to communicate with UE $n$, $\beta$ is the path-loss exponent, $B_m$ is the bandwidth of BS $m$, and $N_0$ is the thermal noise. The long-run average transmission rate is defined as
\begin{equation}
	\bar{L} = \lim\limits_{t \to \infty} \frac{1}{t} \int\limits_0^t  \sum\limits_{n \in \mathscr{N}} O_{m,n}(t)r_{m,n}(t),
	\label{eqn:initial_optimization}
\end{equation}
where $O_{m,n}(t) = 1$ if UE $n$ is associated to BS $m$ at $t$ and $O_{m,n}(t) = 0$ otherwise. It is straightforward to observe that all $O_{m,n}(t)$ values are determined by the association decision $A_i$ taken before $t$. 

Our optimization problem can be formulated as 
\begin{equation}
\begin{aligned}
  \max \quad & \bar{L} \\
	 \text{s.t.:} \quad  \sum\limits_{m \in \mathscr{M}} O_{m,n}(t) = 1, & \quad \forall n \in \mathscr{N}, t > 0; \\
	 O_{m,n}(t) - O_{m,n} &(\text{last}(t)^-) = 0,\\ \forall n \neq U(t), \forall & m \in \mathscr{M},  t > 0; \\
 \sum_{m \notin D(last_{t})}O_{m,n}(t) = 0,& \quad \text{for } n  = U(t), t > 0; \\
 \text{cnt}_m(t) - \sum\limits_{n\in\mathscr{N}} O_{m,n}& (t) = 0, \quad \forall m \in \mathscr{M},
\end{aligned} 
\label{eqn:mobility}
\end{equation}
where $\text{last}(t)^-$ refers to the moment just before $\text{last}(t)$. 

Note that this problem is an NP-hard problem. We provide a sketch of the proof. A sufficient condition to prove that a problem is NP-hard is that a known NP-hard problem can be reduced to the concerned problem~\cite{cormen2009introduction}. In our case, it has been shown in existing research that the static user association optimization problem (the static UA problem), which focuses on associating UEs to BSs in a single moment to maximize network performance or minimize total cost, is NP-hard~\cite{mlika2016user}. In this sense, to show that problem \eqref{eqn:mobility} is NP-hard, we need to show that \eqref{eqn:mobility} is reducible to the static UA problem, such that if a black box solver $S$ can solve the problem~\eqref{eqn:mobility}, $S$ should be able to solve all static UA problems. 

We then consider a virtual scenario to show the connection between problem~\eqref{eqn:mobility} and the static UA problem. Note that this virtual scenario is the only demonstration to show that problem~\eqref{eqn:mobility} is NP-hard and does not apply to the actual implementation of our proposed algorithm. First, we define a virtual BS $v$, which all UEs are associated with at $t = 0$ but will not be available for the association for any UE that has already connected to another BS. We further assume that all UEs are located near the (virtual) boundary of $v$ at $t=0$. Then, we move every UE sequentially until they cross the boundary and trigger a virtual decision epoch. In this sense, the time gap between consecutive virtual decision epochs would be very short, and all UEs would be associated with a non-virtual BS by $t=\delta$, where $\delta$ is an arbitrarily small positive actual number. Suppose our objective is to maximize the transmission rate, the above problem could be formulated either as a static UA (as each UE would only re-associate once) to maximize the instant transmission rate or as an instance of the problem~\eqref{eqn:mobility} to maximize the average transmission rate from $t=0$ to $t=\delta$. Again, as the static UA problem is known to be NP-hard, problem~\eqref{eqn:mobility} is also NP-hard. Due to the space limit, we have proven the NP-hard only by a simplified scenario; yet a complete proof for more general settings can also be readily conducted.


\section{Sequence Q-learning Algorithm}~\label{sec:algorithm}
We now define necessary notations and clarify related concepts to describe the Sequence Q-learning Algorithm (SQA) proposed to solve the problem~\eqref{eqn:mobility} in detail. 

Recall that, for the decision epochs from $t_1$ up to $t_j$, we have a tuple of decisions $A_1, A_2, \cdots, A_j$, with each decision $A_i$ representing the BS associated by UE $U(t_i)$ at decision epoch $t_i$. Define an operation $\oplus$ as attaching an element to the end of a tuple, namely $(x_{1}, x_{2}, ... ,x_{n}) \oplus$ $ x_{n+1} = (x_{1}, x_{2},...,x_{n},x_{n+1})$. 

We further denote $S_{k} = (A_1, A_2, ...,  A_{k-1})$ as the set of decisions that have been made before the decision epoch $t_K$. Our aim can be redefined as identifying the best decision series  $(A_{k}, A_{k+1}, ...,  A_{T}|S_{k})$ at every $t_k$. Therefore, problem~\eqref{eqn:mobility} can be considered as a Markov Decision Process (MDP). For any decision series $S_k$, we further use $C_{i}^{(k)}$ to represent the set of all possible decision series up to the decision epoch $t_i$, where $i > k$. That is, $ C_{i}^{(k)} = \left\lbrace S_{k} \oplus d_{k} \oplus d_{k+1} \oplus ...\oplus d_{i-1} : d_j \in D(t_j) \right\rbrace $. 
For convenience, $C_{i}^{(k)}$ represents an ordered set sorted by the index of $d_k$, then by $d_{k+1}$, and so on. We then denote the $j$th element in $C_{i}^{(k)}$ as $c_{i,j}^{(k)}$. Note that the number of elements in $C_{i}^{(k)}$ is $\sum_{j=1}^{\prod_{l=k}^{i-1} |D_{l}|}$. 

As problem~\eqref{eqn:mobility} can be considered as an MDP, it is straightforward to attempt to solve it by the RL techniques. As one of the classical RL techniques, Q-learning performs well in MDPs with low dimension discrete action space~\cite{watkins1992q}. However, the relatively high-dimensional state space of  problem~\eqref{eqn:mobility} prohibits the direct application of the classical Q-learning method. To address this issue, we propose the SQA based on Q-learning. While the SQA retains the advantages of Q-learning, it is much better in handling high-dimensional state space problems, as we will demonstrate later.

In the classical Q-learning~\cite{dietterich2000hierarchical}, the optimal action-value function for an action $a$, when the current state is $s$, is defined as
\begin{equation}
    \begin{aligned}
     Q^{*}(s,a)=\sum_{s^{'}}P
	\bigg( \left. s^{'}\right|s,& a \bigg) 
	\bigg[R\left( \left.s^{'}\right|s,a \right) \\
	&+ \gamma\cdot \max_{a^{'}} Q^{*}\left( s^{'},a^{'}\right) \bigg], 
	\label{eqn:classicalQ}   
    \end{aligned}
\end{equation}
where $P\left( s^{\prime} \mid s,a \right)$ is the transition probability that the next state is $s^{\prime}$ when taking action $a$ under current state $s$, $R\left(s^{\prime} \mid s,a \right)$ is the reward associated with the action and the transition, and $\gamma \in (0,1)$ is the discount rate. 

Similar to~\eqref{eqn:classicalQ}, we use $Q_{i}^*(S_k, A_i)$ to indicate the optimal action-value function under given state history $S_k$, which can be achieved by connecting UE $U(t_i)$ to BS $A_i$ at decision epoch $t_i$. Specifically, 
\begin{equation}
	\begin{aligned}
		&Q_{i}^* \left(S_{k}, A_{i}\right)
		=E_{\textbf{w}}\left[Q^*_i \left(S_{i}, A_{i}\right)\right] =\sum_{j=1}^{\prod_{l=k}^{i-1} |D_{l}|} w_{i,j}  Q^*_i \left(c_{i,j}^{(k)}, A_{i}\right) \\
		&=\sum_{j=1}^{\prod_{l=k}^{i-1} |D_{l}|} w_{i,j}\left(R_{c_{i,j}^{(k)}, A_{i}}+\gamma \max_{A_{i+1}} Q^*_{i+1} \left(c_{i,j}^{(k)}\oplus A_{i}, A_{i+1}\right)\right) \\
		&=\sum_{j=1}^{\prod_{l=k}^{i-1} |D_{l}|}  w_{i,j} R_{c_{i,j}^{(k)}, A_{i}} + \gamma \max_{A_{i+1}} Q^*_{(i+1)|i} \left(S_{k}, A_{i+1} \mid A_{i}\right)
		\label{eqn:optimal_action_value1}
	\end{aligned}
\end{equation}
where $\textbf{w}$ is the weighed matrix, $w_{i,j} \in \textbf{w}$ is the weight of $c_{i,j}^{(k)}$, $\gamma$ is the discount rate, and 
\begin{equation}
    \begin{aligned}
     Q^*_{(i+1)|i} \bigg(S_{k}, & A_{i+1} \mid A_{i}\bigg)\\
     &=\sum_{j=1}^{\prod_{l=k}^{i-1} |D_{l}|} w_{i,j} Q^*_{i+1} \left(c_{i, j}^{(k)}\oplus A_{i}, A_{i+1}\right). 
	\label{eqn:detailQ}   
    \end{aligned}
\end{equation}

For each decision epoch $i$, if $Q^*_i(S_k, d_{i,j}) > Q^*_i(S_k, d_{i,l})$, we claim that BS candidate $d_{i,j}$ is better than candidate $d_{i,l}$, where $d_{i,j}, d_{i,l} \in D_i$. Thus, if we claim a BS candidate $b$ is good, it shall have a relatively high $Q^*_i(S_k, b)$ among all candidates in $D_i$. Ideally, there exists an appropriate $\textbf{w}$ that makes $Q^*_i(S_k, d_{i,j}) > Q^*_i(S_k, d_{i,l})$ if $d_{i,j}$ is chosen in global optimal solution and $j \neq l$. We can construct the global optimal solution with the previous instructions by connecting the corresponding UE to the best BS candidate at each decision epoch.

Unfortunately, because it is not realistic to enumerate all the possible decision series up to decision epoch $i$, we cannot obtain an ideal and true $Q$ table and $\textbf{w}$. Thus, we use $\hat{Q}$ to approximate the $Q^*$ table and $\hat{\textbf{w}}$ to approximate \textbf{w}. It is straightforward to vertify that 
\begin{equation}
	\begin{aligned}
		\hat{Q_i}&\left(S_{k}, A_{i}\right)
		 =E_{\hat{\textbf{w}}}\left[\hat{Q_i}\left(S_{i}, A_{i}\right)\right] \\
		&=\sum_{j=1}^{\prod_{l=k}^{i-1} |D_{l}|}  \hat{w}_{i,j} R_{c_{i,j}^{(k)}, A_{i}} + \gamma \max_{A_{i+1}} \hat{Q}_{(i+1)|i} \left(S_{k}, A_{i+1} \mid A_{i}\right).
	\end{aligned} 
	\label{eqn:optimal_action_value2}
\end{equation}

Since it is hard to directly give a reasonable $\hat{\textbf{w}}$, while implementing the algorithm to calculate $\hat{Q}$, we let $\hat{\textbf{w}}$ be obtained from exploration strategy instead of assigning fixed values to $\hat{\textbf{w}}$.

To show how to obtain $\hat{\textbf{w}}$ by the exploration strategy, we need to formulate the exploration process. The transition probability $p_{i,j}$ is introduced for describing the exploration strategy, i.e., the probability of connecting corresponding UE to BS candidate $d_{i,j}$ at the decision epoch $t_i\in T$ ($p_{i,j} = 0$ if $d_{i,j}$ is not feasible). Then, we can formulate the exploration process as 
\begin{equation}
\begin{aligned}
	\hat{Q}_{i+1}\left(S_{k}, A_{i+1}\right) 
	&=\sum_{j=1}^{|D_i|} p_{i,j}\hat{Q}_{(i+1)|i}
	   \left(S_{k}, A_{i+1} \mid d_{i,j}\right) \\
	&=\sum_{j=1}^{\prod_{l=k}^{i} |D_{l}|} \hat{w}_{i,j}  \hat{Q}_{i+1}\left(c_{i+1,j}^{(k)}, A_{i+1}\right). 
\end{aligned}  	\label{eqn:optimal_action_value3}
\end{equation}

Referring to the concept of ideal $\textbf{w}$ we mentioned before, we can infer that a good approximation of the ideal $\textbf{w}$ should give high weight to those real good BS candidates at each decision epoch, which means $\hat{Q}_{(i+1)|i}\left(S_{k}, A_{i+1}\right)$ should be dominated by $\hat{Q}_{i+1}\left(c_{i,j}^{(k)}\oplus A_i, A_{i+1}\right)$ of good $A_i$. Thus, in algorithm implementation, we can use the Monte Carlo method to sample good decision series and use their return (the concept in MDP) to estimate $\hat{Q}_{(i+1)|i}\left(S_{k}, A_{i+1}\right)$.

To get good decision series and $\hat{\textbf{w}}$, we need an exploration strategy that assigns a high transition probability to good BS candidates at each decision epoch. Specifically, the transition probability $p_{i,j}$ is defined as
\begin{equation}
	p_{i,j} = \frac{\epsilon^{\phi_i(j)}}{\sum_{l=1}^ {|D_i|}\epsilon^{\phi_i(l)}}, \label{eqn:transition_pro}
\end{equation}
where $\epsilon$ is a preset hyper parameter, and $\phi_i(j)$ is the number of actions in $D_i$ whose $\hat{Q}$ value is less than $\hat{Q}_i\left(S_k, d_{i, j}\right)$. Then, we can design the update rule as
\begin{equation}
	\begin{aligned}
		\hat{Q_i}\left(S_{k}, A_{i}\right)
		\leftarrow 
		&\hat{Q_i}\left(S_{k}, A_{i}\right) 
		+\alpha \gamma \max _{A_{i+1}} \hat{Q}_{(i+1)|i}\left(S_{k}, A_{i+1}\right)\\
		&-\alpha \hat{Q_i}\left(S_{k},A_{i}\right)
		+\text {R}_{S_i} ,	
	\end{aligned} \label{eqn:optimal_action_value4}
\end{equation}
where $\alpha$ is the learning rate. $\text {R}_{S_i}$ is the sum transmission rate between $t_i$ to $t_{i+1}$. Assume that decision series $s_i \in c_{i}^{(k)}$ are made the same, then 
\begin{equation}
	\text {R}_{s_i} 
	= \sum_{l=1}^{N}
	\displaystyle{\int_{t_{i+1}-t_i}^{t_{i+1}-t_i}}
	r(t)\Delta t.
	\label{sum_transmission_rate}
\end{equation}

The pseudo-code of the SQA is presented in Algorithm~\ref{Algorithm1}. In Algorithm~\ref{Algorithm1}, we set a parameter $STEP$ to indicate the maximum number of steps for exploration from any decision epoch, which could control the running time of the SQA. By adjusting the values of $STEP$ and the discount rate $\gamma$, we may fine-tune the performance and running time according to specific network scenarios and UE trajectories.
\renewcommand{\algorithmicrequire}{\textbf{Input:}}  
\renewcommand{\algorithmicensure}{\textbf{Output:}} 
\begin{algorithm}[h]
	\caption{Sequence Q-Learning} 
	\begin{algorithmic}[1]
		\Require
		$\hat{Q_i}\left(S_{k-1}, d_{i,j}\right),\ D_i,\  S_{k-1}$
		\Statex $\qquad \qquad\qquad \ \quad \quad 
		i \in \left\lbrace k,k+1,...,T\right\rbrace, d_{i,j} \in D_i  $
		\Ensure
		$\hat{Q_i}\left(S_{k},d_{i,j}\right) \quad i \in \left\lbrace k,k+1,...,T\right\rbrace $
		\For{$i=k$ to $T$}
			\State $\hat{Q_i}\left(S_{k}, d_{i,j}\right) 
			 \longleftarrow 
			 \hat{Q_i}\left(S_{k-1}, d_{i,j}\right) $
		\EndFor
		\For{$j=0$ to max-iteration}
		    \State $step \longleftarrow 0$
			\State Explore $(S_i)$
			\State update $p_{i,j}$ by 
			        $\hat{Q_i}\left(S_{k}, d_{i,j}\right)$ 
	    \EndFor
	    \State \quad
	    \State \textbf{Procedure}  Explore $(S_i: \text{decision series up}$ 
	    \Statex $\text{ to decision epoch}\ i )$
	    \If{$\left( i==T\right) $ or $\left( step > STEP \right)$}
 			\State 
	    	\Return 0
	    \EndIf
	    \State choose BS $A_i$ by $p_i$
	    \State $step \longleftarrow step + 1 $ 
	    \State  $S_{i+1}\longleftarrow S_i \oplus A_i $
	    \State  $\text{return}_{\text{MDP}} \longleftarrow \text{R}_{S_{i+1}} + \gamma \text{Explore}(S_{i+1}) $
	    \State  $\hat{Q}_{(i+1)|i}\left(S_{k}, A_{i+1}\right) \longleftarrow \hat{Q}_{(i+1)|i}\left(S_{k}, A_{i+1}\right)$
	    \Statex $\qquad \qquad \qquad +\alpha \left[\text{return}_{\text{MDP}}-\hat{Q}_{(i+1)|i}\left(S_{k}, A_{i+1}\right) \right] $
	    \State  $\hat{Q}_{i}\left(S_{k}, A_{i}\right) \longleftarrow
	    	\hat{Q}_{i}\left(S_{k}, A_{i}\right)
	    	+\alpha \bigg[\text{R}_{S_{i+1}}-\hat{Q}_{i}\left(S_{k}, A_{i}\right)$
	    \Statex$\qquad \qquad\qquad+ \gamma \max_{A_{i+1}} \hat{Q}_{(i+1)|i}\left(S_{k}, A_{i+1}|A_i\right)\bigg] $
	    \State 
	    \Return $\text{R}_{S_{i+1}}+\gamma \text{return}_{\text{MDP}}$
	\end{algorithmic} \label{Algorithm1}
\end{algorithm}

\begin{figure}[htp!]
	\centering
	\includegraphics[width=0.35\textwidth]{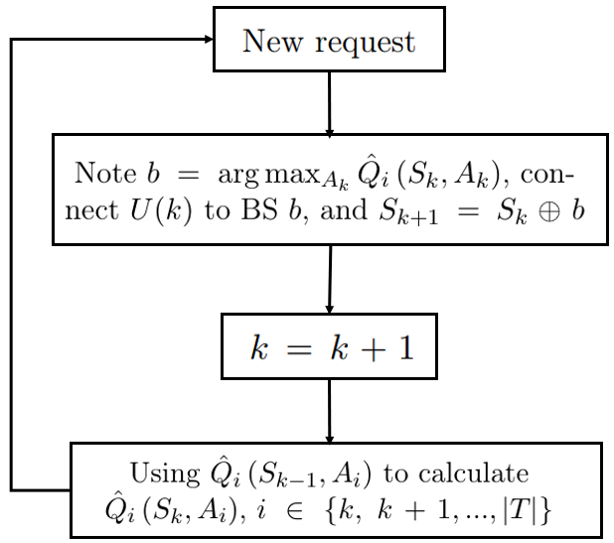}
	\caption{Flowchart of applying Algorithm 1 to solve our problem.}
	\label{fig:Flowchart of using the Algorithm 1 into our scenario}
\end{figure}

Next, we illustrate how we apply the SQA to solve \eqref{eqn:mobility}, with the flowchart shown in Fig.~\ref{fig:Flowchart of using the Algorithm 1 into our scenario}. 
It is worth noting that SQA is called between two decision epochs instead of after each new request. In this way, the association decision is made immediately after each request without any delay because the decision can be made by only checking the $\hat{Q}$ table.

Finally, a brief explanation of the fact that SQA has a polynomial time complexity is shown in the following. For each decision, SQA needs to run $I$ iterations, and the number of steps in each iteration is $STEP$. In each step, a linear time complexity $O(|\mathscr{M}|)$ is required to determine the next action while exploring and updating corresponding $\hat{Q}$ values once using backtracking. Therefore, the total time complexity of SQA for one decision is $O(I \cdot STEP \cdot |\mathscr{M}|)$, which is polynomial.

\section{Numerical Experiments} \label{sec:experiment}
\subsection{Experiment setup}
\subsubsection{Network map and UE trajectory generation}

We consider a square network map with $Z\times Z$ grids, with grid lines representing roads, the length of which is set to $200$m. BSs are deployed at each intersection of grid lines. We consider an area of $1600$m $\times 1600$m, which corresponds to  $8 \times 8$ grids with $64$ BSs. Each BS covers a circular area with a radius of $300$m. A rectangular building with random length and width (not exceeding the road length) is randomly placed at each grid. 

Due to the nature of mmWave communications, a building will block any transmissions between BSs and UEs if the building is on the Line-of-Sight (LoS) transmission path. The bandwidth $B_m$ of all BSs is uniformly set to 10MHz, the average transmit power of all BS-UE pairs is also uniformly set as $P_{m,n}= P = 30$dBm. We also set $N_0 = -90$dBm, and $\beta = 3$. 

We further assume that UEs are only distributed on the roads and only move along them. UEs move at an average speed of $15$m/s towards a random possible direction, which ultimately form the predetermined trajectories. We simulate SQA and the benchmarks in four scenarios with different densities of UEs. 

\subsubsection{Benchmark approaches} 
We compare the performance of SQA with four state-of-the-art benchmarks: SNR-based handover (SBH), Rate-based handover (RBH), Learning-based handover (LBH)~\cite{khosravi2020learning}, and SMART~\cite{sun2017smart}. SBH and RBH are greedy algorithms which always choose the BS with the highest SNR or instant transmission rate at a decision epoch, respectively. SMART and LBH are two recently proposed RL approaches. SMART aims to 
control the frequency of handovers, subject to the amount of incremental total transmission rate after a handover. LBH uses Q-learning to optimize the long-term transmission rate, assuming  that a constant bandwidth is offered to every UE.

\subsubsection{Experiment initialization and parameter settings} 
Note that an initial $\hat{Q}$ table is necessary to start the SQA. We use a greedy algorithm to initialize the $\hat{Q}$ table. In addition, we set $\epsilon = 3$, $\alpha = 0.01$, $\gamma = 1.0$, $STEP = \frac{N}{2}$ (half of the total number of UEs), and the number of iterations is set to 100. 

\subsection{Numerical results}
\begin{table}
	\caption{Performance comparison of 5 algorithms in 4 scenarios with different UE densities.}
	\centering
	\begin{tabular}{|c|c|c|c|c|c|}
		\hline
		\multicolumn{6}{|c|}{Very low density of UEs (BS:UE=64:512)} \\
		\hline
		& SBH   & RBH   & LBH   & SMART & SQA \\
		\hline
		$\bar{L}$ (10Mbit/s) & 77.60 & 88.09 & 81.64 & 84.67 & 96.26 \\
		\hline
		$X_n$ (s) & 14.42 & 13.40 & 17.96 & 14.89 & 14.76 \\
		\hline
		\multicolumn{6}{|c|}{Moderately low density of UEs (BS:UE=64:1024)} \\
		\hline
		& SBH   & RBH   & LBH   & SMART & SQA \\
		\hline
		$\bar{L}$ (10Mbit/s) & 81.33 & 85.80 & 81.71 & 84.62 & 93.67 \\
		\hline
		$X_n$ (s) & 14.50 & 13.57 & 17.96 & 14.62 & 14.59 \\
		\hline
		\multicolumn{6}{|c|}{Moderately high density of UEs (BS:UE=64:1536)} \\
		\hline
		& SBH   & RBH   & LBH   & SMART & SQA \\
		\hline
		$\bar{L}$ (10Mbit/s) & 80.51 & 85.81 & 82.32 & 83.56 & 87.28 \\
		\hline
		$X_n$ (s) & 14.39 & 13.28 & 17.79 & 14.14 & 14.37 \\
		\hline
		\multicolumn{6}{|c|}{Very high density of UEs (BS:UE=64:2048)} \\
		\hline
		& SBH   & RBH   & LBH   & SMART & SQA \\
		\hline
		$\bar{L}$ (10Mbit/s) & 82.07 & 84.29 & 84.28 & 84.01 & 85.76 \\
		\hline
		$X_n$ (s) & 14.45 & 14.04 & 17.90 & 14.52 & 14.50 \\
		\hline
	\end{tabular}%
	\label{tab1}%
\end{table}%

\begin{figure*}[htp!]
\centering
\subfigure[Very low density of UEs (BS:UE=1:8)]{
	\includegraphics[width=6.5cm]{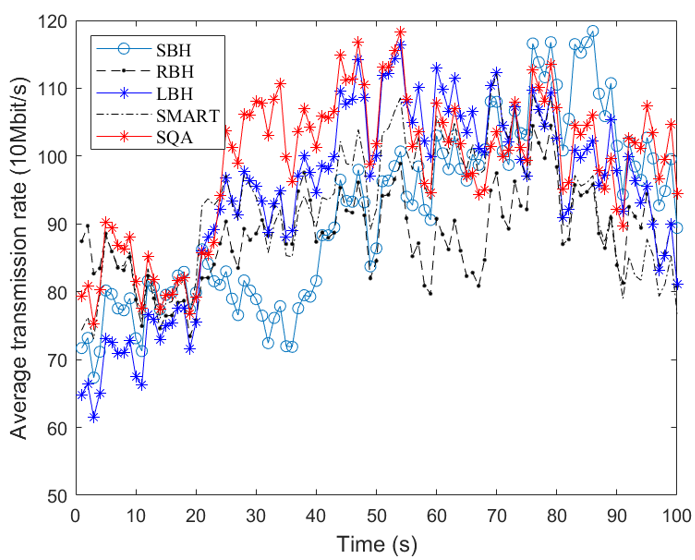}%
}
\quad 
\subfigure[Moderately low density of UEs (BS:UE=1:16)]{
	\includegraphics[width=6.5cm]{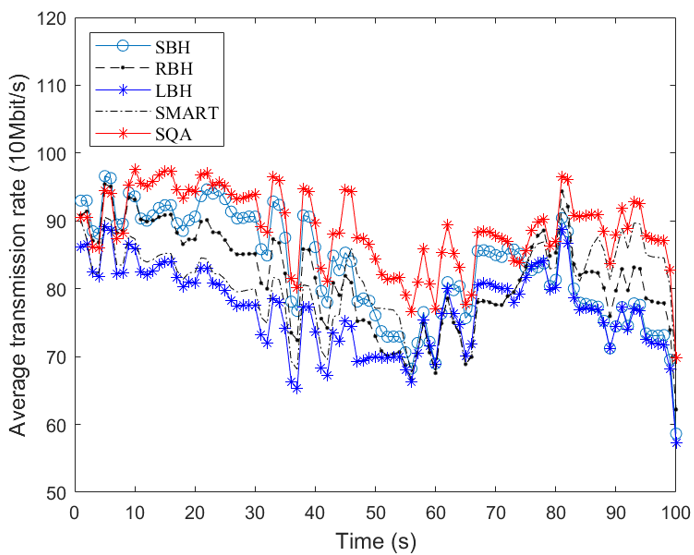}%
}
\quad
\subfigure[Moderately high density of UEs (BS:UE=1:24)]{
	\includegraphics[width=6.5cm]{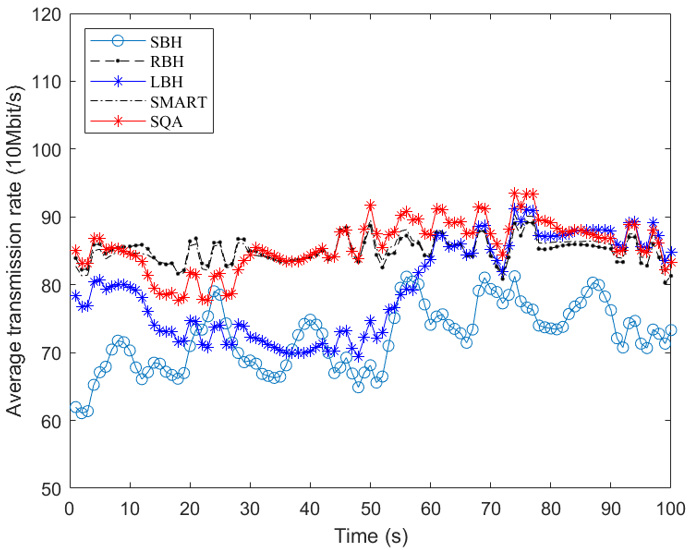}%
}
\quad
\subfigure[Very high density of UEs (BS:UE=1:32)]{
	\includegraphics[width=6.5cm]{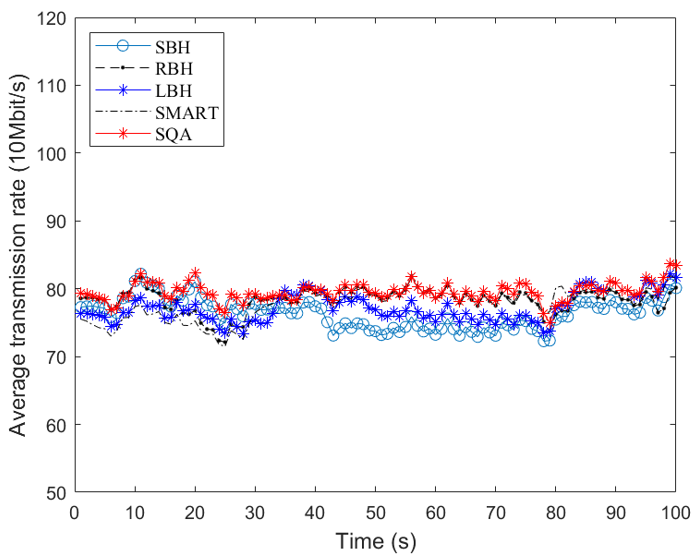}%
}
\caption{The comparison of instant transmission rate achieved by different algorithms in the first 100 seconds.}
\label{fig:Comp}
\end{figure*}%
Table~\ref{tab1} shows the simulation results. $\bar{L}$ is the long-run average transmission rate of the network, and $X_n$ is the average time between two successive handovers for a single UE. Recall that the triggering condition for handover events is SNR, which we assumed to be influenced only by the distance between the UE and the BS. Therefore, $X_n$ under a particular algorithm is not affected by the densities of UEs, but determined by the average speed of UEs.   

In terms of $\bar{L}$, SQA outperforms four  benchmarks in the four scenarios. Notably, in low density scenarios, SQA can improve over $10\%$ throughput than other algorithms. For $X_n$, the performances of SQA and two RL-based benchmarks are close. We also observe that, by utilizing the future trajectory information of UEs, SQA and LBH attain significantly better $\bar{L}$ and $X_n$ than the other three algorithms, respectively.

The greedy algorithms fail to attain the balance between $\bar{L}$ and $X_n$. While RBH can also achieve comparatively  high $\bar{L}$, it incurs the most frequent handovers as reflected in the lowest $X_n$ among all algorithms. On the other hand, LBH leads to the least frequent handovers, but its $\bar{L}$ is not satisfactory under scenarios with relatively low UE densities. 

To demonstrate the results more intuitively, we show the instant average transmission rate achieved by different algorithms for the first 100 seconds in Fig.~\ref{fig:Comp}. It can be observed that when the UE density is low, SQA can keep the instant rate at a relatively high level and is remarkably superior to other algorithms. While the differences between algorithms are less significant as the UE density increases, SQA still outperforms other algorithms most times. Combining the long-run average results demonstrated in Table I, SQA is able to achieve remarkable improvements for both short and long periods of time. 

\section{Conclusion}~\label{sec:conclusion}
This paper investigated the problem of optimizing the long-run average transmission rate with a feasible user association scheme in mmWave networks while controlling the number of handovers to a reasonable level. To solve this NP-hard problem, we proposed the SQA, which has a polynomial time complexity concerning the number of BSs and UEs. Simulation results demonstrated that the SQA outperforms all benchmark algorithms, including two greedy approaches and two recently proposed RL-based algorithms, in a range of scenarios with different densities of UEs.

\section*{Acknowledgement}

The work described in this paper is partly supported by College Research Grant from BNU-HKBU United International College R201911, and Zhuhai Basic and Applied Basic Research Foundation Grant ZH22017003200018PWC. The authors sincerely thank 4 anonymous reviewers for their constructive comments and suggestions that helped to improve the paper.

\bibliographystyle{IEEEtran}
\bibliography{IEEEabrv, mybibliography} 
\end{document}